\documentclass{amsart}

\usepackage{amsmath,amsthm,amssymb, amsaddr}
\usepackage{graphicx}
\usepackage{url}

\usepackage{setspace}

\theoremstyle{plain}
\newtheorem{theorem}{Theorem}[section]

\theoremstyle{definition}

\theoremstyle{remark}

\newcommand{\real}{\mathbb{R}}
\newcommand{\cg}[1]{$\mathrm{cG}(#1)$}
\newcommand{\dg}[1]{$\mathrm{dG}(#1)$}
\newcommand{\dt}{\,\mathrm{d}t}
\newcommand{\ds}{\,\mathrm{d}s}
\newcommand{\inner}[2]{\langle #1, #2 \rangle}
\newcommand{\RR}{\bar{R}}
\newcommand{\emach}{\epsilon_{\textrm{mach}}}
\newcommand{\nmach}{n_{\textrm{mach}}}
\newcommand{\EE}{\mathbf{E}}
\newcommand{\ED}{\mathbf{E}_D}
\newcommand{\EG}{\mathbf{E}_G}
\newcommand{\EC}{\mathbf{E}_C}

\renewcommand{\k}{\Delta{}t}

\let\oldmarginpar\marginpar
\renewcommand\marginpar[1]{\-\oldmarginpar[\raggedleft\footnotesize #1]%
{\raggedright\footnotesize #1}}

\title{Quantifying the computability of the Lorenz system}
\author{Benjamin Kehlet and Anders Logg}
\address{Center for Biomedical Computing, Simula Research Laboratory \\ P.O.Box 134, 1325 Lysaker, Norway}


\begin{document}

\begin{abstract}
  It is well known that the computation of accurate trajectories of the
  Lorenz system is a difficult problem. Computed solutions are very
  sensitive to the discretization error determined by the time step size
  and polynomial order of the method, as well as round-off errors.

  In this work, we show how round-off errors limit the computability
  of the Lorenz system and quantify exactly the length of intervals
  over which solutions can be computed, expressed in terms of the
  floating point precision. Using adjoint-based \emph{a~posteriori}
  error analysis techniques, we estimate the stability of computations
  with respect to initial data, discretization, and round-off errors,
  respectively.

  The analysis is verified by computing an accurate solution on the
  time interval $[0, 1000]$ using a very high order (order 200) finite
  element method and very high floating point precision (400 digits).
\end{abstract}

\maketitle

\section{Introduction}

In a classic paper from 1963~\cite{lorenz_deterministic_1963}, Edward
Lorenz studied the computability of a simple system of three ordinary
differential equations,
\begin{equation} \label{eq:lorenz}
  \left\{
   \begin{aligned}
    \dot{x} &= \sigma (y-x), \\
    \dot{y} &= rx - y -xz,   \\
    \dot{z} &= xy - bz,      \\
   \end{aligned}
  \right.
\end{equation}
where $\sigma = 10$, $b = 8/3$, and $r = 28$. Lorenz computed
numerical solutions of the system~\eqref{eq:lorenz} and found the
solutions to be very sensitive to changes in initial data. The
equations had been devised by Lorenz as a simple model of atmospheric
flow, based on a truncated Fourier expansion of the partial
differential equations governing Rayleigh--B\'enard
convection~\cite{rayleigh_convective_1916,lorenz_maximum_1960,saltzman_finite_1962}.
In his paper, Lorenz computed solutions on the interval $[0, 60]$. As
we shall see below, the Lorenz system is not computable on the
equipment that was available to Lorenz in~1963 beyond time $T \approx
25$.

It is known that, given enough resources, the Lorenz system is
computable over arbitrarily long time intervals. However, one may
easily (and falsely) come to the conclusion that the Lorenz system is
computable only over very short time intervals, either by numerical
experiments or by a simplistic analysis. Indeed, a standard
\emph{a~priori} error estimate indicates that the growth rate of the
error is
\begin{equation} \label{eq:apriori}
  \|e(T)\| \leq C e^{LT} \epsilon,
\end{equation}
where $\|e(T)\|$ denotes some norm of the error at the final time~$T$,
$L$ is the Lipschitz constant of~(\ref{eq:lorenz}), and $\epsilon$ is
the size of the residual or local truncation error in a numerical
solution of~\eqref{eq:lorenz}. The Lipschitz constant is of size $L
\approx 33$ which indicates that solutions are not computable beyond
$T \approx 1.1$, even if the residual is close to machine precision
($\emach \sim 10^{-16}$ on most computers).\footnote{The value of the
  Lipschitz constant was computed as the maximum $l^2$-norm of the
  Jacobian~$J = \partial f / \partial u$ of the right-hand side~$f$ of
  the Lorenz system over the interval~$[0, 1000]$.} However, the
estimate~\eqref{eq:apriori} is overly pessimistic; it is well known
that solutions of the Lorenz system may be computed on short time
intervals. In fact, one may easily compute accurate solutions over
time intervals of length $T = 25$ with any standard ODE solver.

In~\cite{estep_pointwise_1998}, it was demonstrated that the Lorenz
system is indeed computable on intervals of moderate length ($T = 30$)
on a standard desktop computer. The computability of the Lorenz system
was linked to the growth of a \emph{stability factor} in an \emph{a
  posteriori} estimate of the error at the final time. It was shown
that the growth rate of the stability factor is non-constant. On
average the growth is exponential but with a rate much smaller than
indicated by~\eqref{eq:apriori}.

In~\cite{logg_multi-adaptive_2003}, the computability of the Lorenz
system was further extended to $T = 48$ using high order ($\|e(T)\|
\sim \k^{30}$) finite element methods. As we shall see below, this is
the ``theoretical limit'' for computations with 16~digit
precision. Solutions over longer time intervals have been computed
based on shadowing (the existence of a nearby exact solution),
see~\cite{coomes_rigorous_1995}, but for unknown initial data. Other
related work on high-precision numerical methods applied to the Lorenz
system include~\cite{viswanath2004fractal}
and~\cite{jorba_software_2005}. For an overview of some recent results
obtained with high-precision numerical methods, we also refer
to~\cite{Bailey201210106}.

In this paper, we study and quantify the computability of the
Lorenz system. In particular we answer the following fundamental
question: \emph{How far is the Lorenz system computable for a given
  machine precision?}

As we shall see, obtaining a sequence of converging approximations for
the solution of the Lorenz system is non-trivial. In particular, such
a sequence of solutions cannot be obtained by simply decreasing the
size of the time step; see for example~\cite{teixeira2007time}. This
has led to misconceptions regarding the computability of the Lorenz
system; see for example~\cite{yao_comment_2008}. To obtain a sequence
of converging solutions, one must also control the effect of round-off
errors. This was also noted by Lorenz~\cite{lorenz_reply_2008} in a
response to~\cite{yao_comment_2008}.

In this manuscript, we define computability as the length $T$ of the
maximum time interval $[0, T]$ on which a solution is computable to
within a given precision $\epsilon > 0$ using a given machine
precision $0 < \emach < \epsilon$; that is, the maximum $T$ such that
$\inf_U \|u - U\|_{L^{\infty}(0, T; l^{\infty})} \leq \epsilon$, where
the infimum is taken over all numerical approximations~$U$ of the
exact solution~$u$ computed with some time-stepping method and machine
precision $\emach$ (as made more precise in
Section~\ref{sec:erroranalysis}). If the computability $T_{\epsilon} =
T_{\epsilon}(\emach)$ does not depend strongly on $\epsilon$, we write
$T = T(\emach)$. As we shall see, this is the case for the Lorenz
system as a result of exponential growth of errors as function of the
final time $T$. The definition of computability $T(\emach)$ is closely
related to the definition of a critical predictable time $T_c$ in
\cite{liao2009reliability} and the definition of a decoupling time
$\hat{T}$ in \cite{teixeira2007time}.

\section{Numerical method and implementation}
\label{sec:methods}

We consider the numerical solution of general initial value problems
for systems of ordinary differential equations,
\begin{equation} \label{eq:u'=f}
  \begin{split}
    \dot{u}(t) &= f(u(t),t), \quad t \in (0, T],
      \\ u(0) &= u_0.
  \end{split}
\end{equation}
The right-hand side $f : \real^N \times [0, T] \rightarrow \real^N$ is
assumed to be Lipschitz continuous in~$u$ and continuous in~$t$. Our
objective is to analyze the error in an approximate solution~$U : [0,
  T] \rightarrow \real^N$, for example a numerical solution of the
Lorenz system.

The continuous and discontinuous Galerkin methods \cg{q} and \dg{q}
are formulated by requiring that the residual $R = \dot{U} -
f(U,\cdot)$ be orthogonal to a suitable space of test functions. By
making a piecewise polynomial Ansatz, the solution may be computed on
a sequence of intervals partitioning the computational domain $[0, T]$
by solving a system of equations for the degrees of freedom on each
consecutive interval. For a particular choice of numerical quadrature
and degree~$q$, the \cg{q} and \dg{q} methods both reduce to standard
implicit Runge--Kutta methods.

In the case of the \cg{q} method, the numerical solution~$U$ is a
continuous piecewise polynomial of degree~$q$ that on each
interval~$[t_{n-1}, t_n]$ satisfies $\int_{t_{n-1}}^{t_n} R \, v \dt =
  0$ for all $v \in \mathcal{P}^{q-1}([t_{n-1}, t_n])$.

The results were obtained using the finite element
package~DOLFIN~\cite{logg_dolfin:_2009,fenics:book} version~0.9.2 together with
the multi-precision library GMP~\cite{gmp}. For a detailed discussion
on the implementation, we refer to~\cite{kehlet_analysis_2010}. The
source code as well as scripts to reproduce all results presented in
this manuscript are available on request.

\section{Error analysis}
\label{sec:erroranalysis}
The error analysis is based on the solution of an auxiliary \emph{dual
  problem}. The dual (adjoint) problem takes the form of an initial
  value problem for a system of linear ordinary differential
  equations,
\begin{equation} \label{eq:dual}
  \begin{split}
    - \dot{z}(t) &= \bar{A}^{\top}(t) z(t), \quad t \in [0, T), \\
    z(T) &= z_T.
  \end{split}
\end{equation}
Here, $\bar{A}(t) = \int_0^1 \frac{\partial f}{\partial u}(sU(t) +
(1-s)u(t), t) \ds$ denotes the Jacobian matrix of the right-hand
side~$f$ averaged over the approximate solution~$U$ and the exact
solution~$u$.

The Lorenz system is quadratic in the primal variable~$u$. Hence, the
average in $\bar{A}$ corresponds to evaluating the Jacobian matrix at
the midpoint between the two vectors $U(t)$ and $u(t)$. It follows
that the dual problem of the Lorenz system is
\begin{equation} \label{eq:lorenz,dual}
  \left\{
  \begin{aligned}
    -\dot{\xi}   &= -\sigma\xi + (r-\bar{z})\eta + \bar{y}\zeta, \\
    -\dot{\eta}  &= \sigma\xi - \eta + \bar{x}\zeta, \\
    -\dot{\zeta} &= -\bar{x}\eta - b\zeta, \\
  \end{aligned}
  \right.
\end{equation}
where $z = (\xi, \eta, \zeta)$ denotes the dual solution and
$(\bar{x}, \bar{y}, \bar{z}) = (U + u)/2$.

In \cite{kehlet_aposteriori_2013}, we prove the following
\emph{a~posteriori} error estimate:
\begin{theorem}[Error estimate] \label{th:errorestimate}

  Let $u : [0, T] \rightarrow \real^N$ be the exact solution
  of~\eqref{eq:u'=f} \emph{(}assuming it exists\emph{)}, let $z : [0,
    T] \rightarrow \real^N$ be the solution of~\eqref{eq:dual}, and let $U
  : [0, T] \rightarrow \real^N$ be any piecewise smooth approximation of
  $u$ on a partition $0 = t_0 < t_1 < \cdots < t_M = T$ of $[0, T]$,
  that is, $U \vert_{(t_{m-1}, t_m]} \in
    \mathcal{C}^{\infty}((t_{m-1}, t_m])$ for $m=1,2,\ldots,M$ \emph{(}$U$ is
      left-continuous\emph{)}.

  Then, for any $p \geq 0$, the following error estimate holds:
  \begin{equation*}
    \inner{z_T}{U(T) - u(T)} = \ED + \EG + \EC,
  \end{equation*}
  where
  \begin{equation*}
    \begin{split}
      |\ED| &\leq S_D \, \|U(0) - u(0)\|, \\
      |\EG| &\leq S_G \, C_p \max_{[0,T]} \left\{\k^{p+1} (\|[U]\|/\k + \|R\|)\right\}, \\
      |\EC| &\leq S_C \, C_p' \max_{[0,T]} \|\k^{-1} \RR\|,
    \end{split}
  \end{equation*}
  where $C_p$ and $C_p'$ are constants depending only on~$p$.
  The stability factors $S_D$, $S_G$, and $S_C$ are defined by
  \begin{equation*}
      S_D = \|z(0)\|, \quad
      S_G = \int_0^T \|z^{(p+1)}\| \dt, \quad
      S_C = \int_0^T \|\pi z\| \dt.
  \end{equation*}
\end{theorem}

Furthermore the following bound for the computational error is proved in \cite{kehlet_aposteriori_2013}:

\begin{theorem} \label{th:random}
  Assume that the round-off error is a random variable of size
  $\pm\emach$ with equal probabilities. Then the root-mean squared
  expected computational error $\EC$ of Theorem~\ref{th:errorestimate}
  is bounded by
  \begin{equation*}
    (E[\EC^2])^{1/2}
    \leq
    S_{{C_2}} \, \sqrt{C_p'} \frac{\emach}{\min_{[0,T]}\sqrt{\k}},
  \end{equation*}
  where $S_{C_2} = \left(\int_0^T \|\pi z\|^2 \dt\right)^{1/2}$ and
  $C_p'$ is a constant depending only on~$p$.
\end{theorem}
We note that the computational error (accumulated round-off error) is
inversely proportional to (the square root of) the time step; that is,
a \emph{smaller} time step yields a \emph{larger} accumulated
round-off error.

\begin{figure}[htbp]
  \begin{center}
    \includegraphics[width=\textwidth]{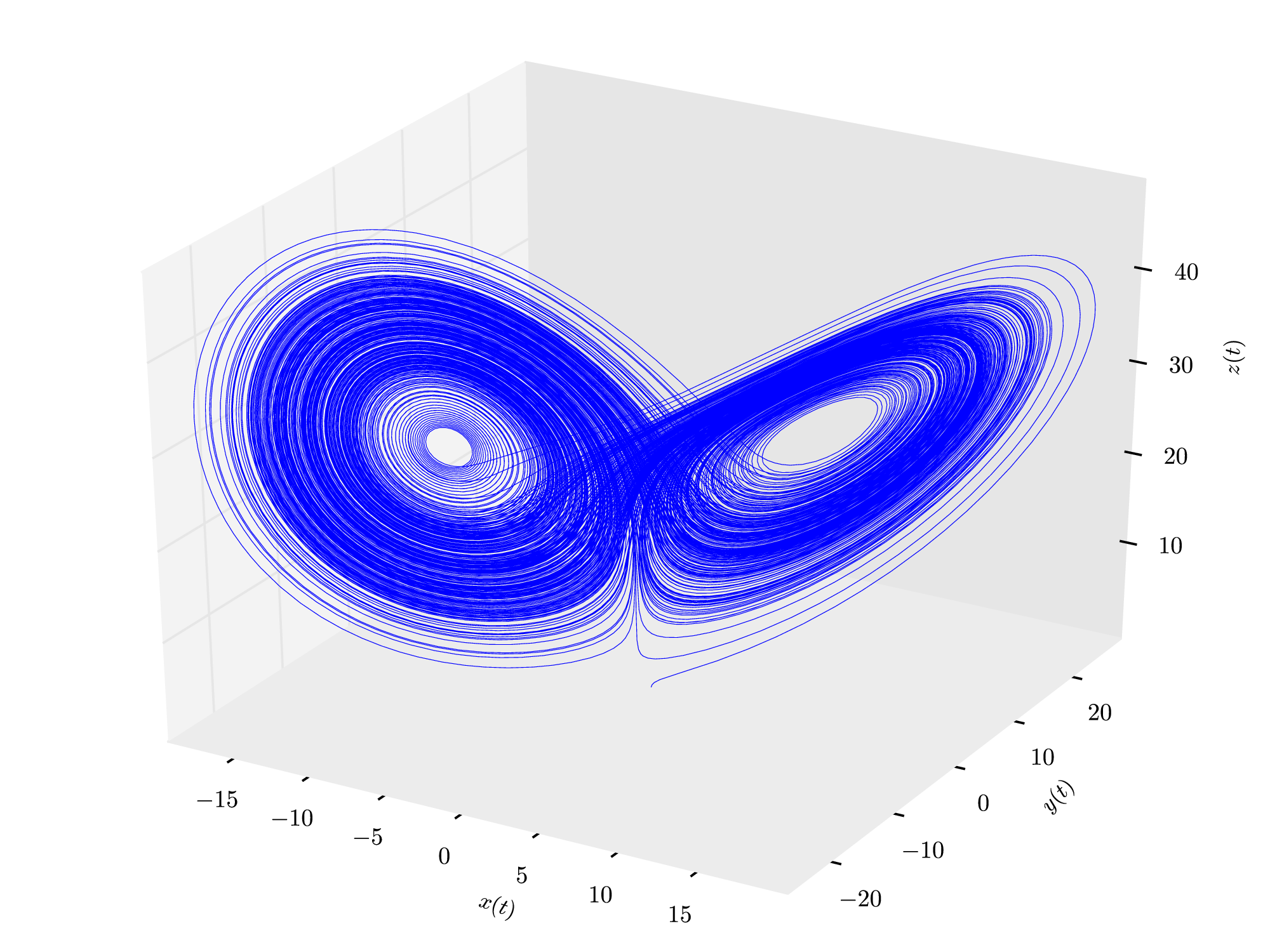}
    \caption{Phase portrait of the solution of the Lorenz system on
      the time interval~$[0, 1000]$ for $u(0) = (1,0,0)$.}
    \label{fig:phaseportrait}
  \end{center}
\end{figure}

\section{Numerical results}
\label{sec:numericalresults}

In this section, we present numerical results in support of
Theorem~\ref{th:errorestimate} and Theorem~\ref{th:random}.

\subsection{Solution of the Lorenz system}

The phase portrait of the solution of the Lorenz system is plotted in
Figure~\ref{fig:phaseportrait}. The solution was computed with
\cg{100}, which is a method of order $2q = 200$, a time step of size
$\k = 0.0037$, 420-digit precision arithmetic\footnote{The requested
  precision from GMP was 420 digits. The actual precision is somewhat
  higher depending on the number of significant bits chosen by GMP.},
and a tolerance for the discrete residual of size $\emach \approx 2.26
\cdot 10^{-424}$. The solution trajectory revolves around one of the
two unstable fixed points $P_{\pm} = (\pm 6\sqrt{2}, \pm 6\sqrt{2},
27)$ for a while and then, seemingly at random, jumps to the other
fixed point. Phase portraits (``attractors'') resembling the phase
portrait of Figure~\ref{fig:phaseportrait} are commonly displayed in
most books on dynamical systems and chaos theory. However, in one way
the phase portrait of Figure~\ref{fig:phaseportrait} is significantly
different. It is the phase portrait of a well-defined dynamical
system, namely the Lorenz system~\eqref{eq:lorenz} with initial
condition $(1, 0, 0)$, not the result of an unspecified discrete map
which includes both the effect of a particular time-stepping scheme
and the unknown effect of round-off errors.

\begin{figure}
  \begin{center}
    \includegraphics[width=\textwidth]{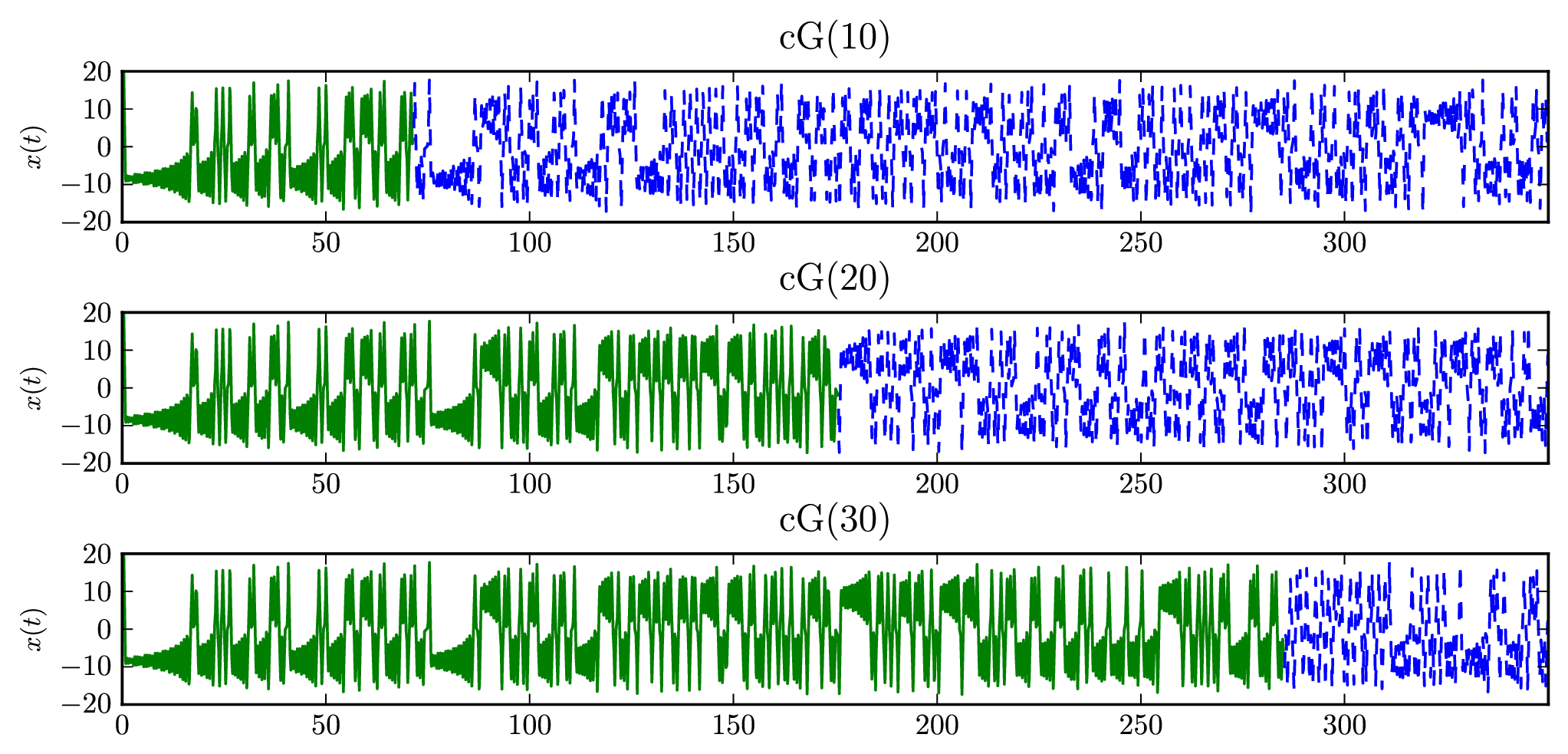}
    \includegraphics[width=\textwidth]{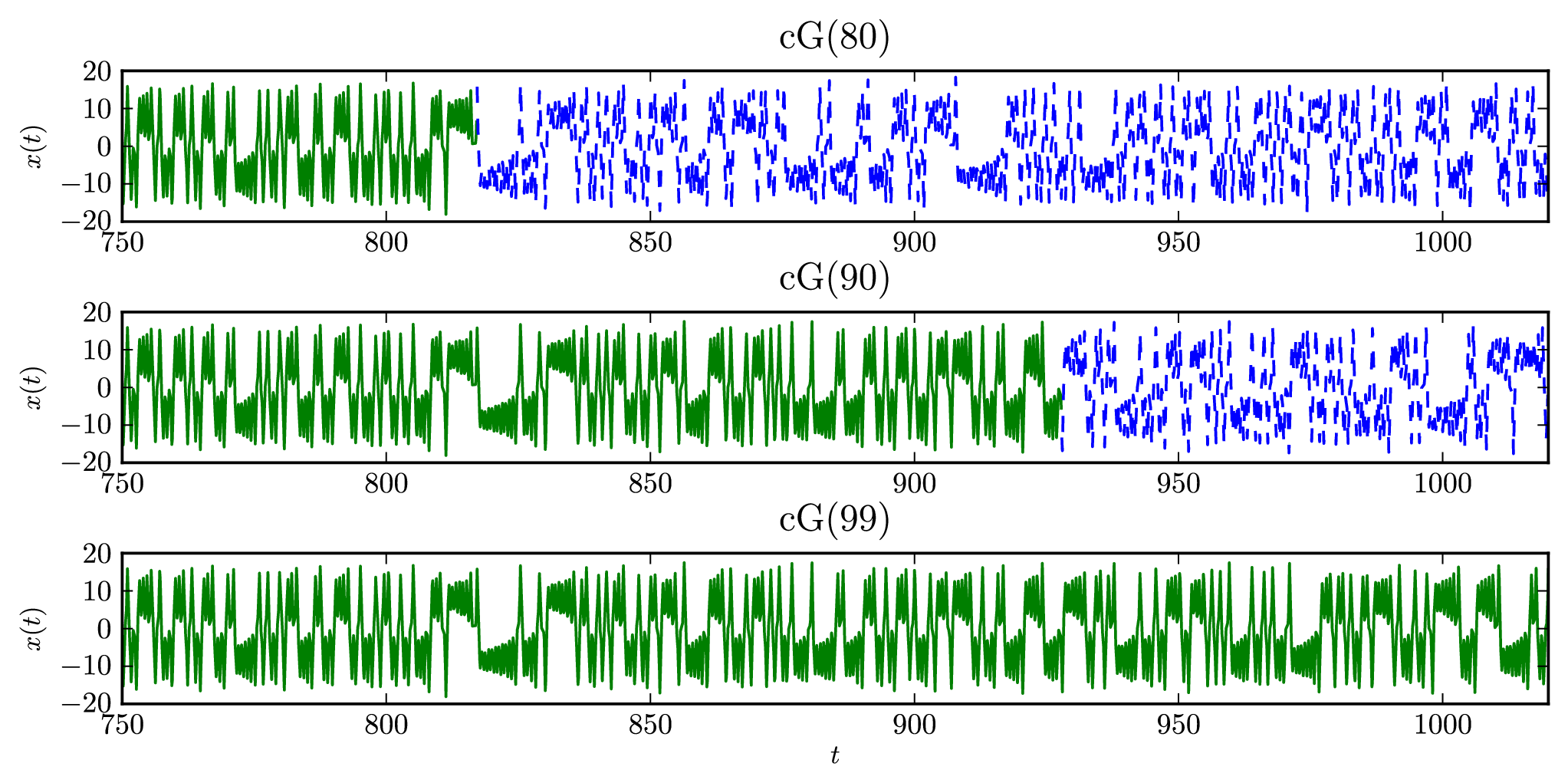}
    \caption{Computed numerical solutions ($x$-component) for the
      Lorenz system with methods of increasing order, starting at
      \cg{10} (a method of order~$20$) and increasing up to \cg{99} (a
      method of order~$198$).}
    \label{fig:convergence}
  \end{center}
\end{figure}

To verify the computed solution, we perform a simple experiment where
we compute the solution with methods of increasing order. The time
step is fixed ($\k = 0.0037$) and so is the arithmetic precision
($420$ digits). By Theorem~\ref{th:errorestimate}, we expect the
discretization error~$\EG$ to decrease exponentially with increasing
order while the computational error~$\EC$ remains bounded. The error
should therefore decrease, until $\EG < \EC$. Since no analytic
solution or other reference solution is available, we compare the
\cg{10} solution with the \cg{20} solution and conclude that when the
two solutions no longer agree to within some tolerance (here
$10^{-16}$), the \cg{10} solution is no longer accurate. The same
experiment is repeated for \cg{20/30}, \cg{30/40}, \ldots,
\cg{90/100}, \cg{99/100}. The solutions are displayed in
Figure~\ref{fig:convergence}. The results indicate that the \cg{99}
solution is accurate on the time interval $[0, 1025]$. Alone, this
does not prove that the \cg{99} is accurate at time $T =
1025$. However, together with the error estimate of
Theorem~\ref{th:errorestimate} and the numerically computed values of
the stability factors presented below, there is strong evidence that
the solution is accurate over $[0, 1025]$.

We emphasize that similar results may be obtained with other numerical
methods and other software. In particular,
Theorem~\ref{th:errorestimate} shows that the solution is computable
with any solver that (i) discretizes the equations with high order and
(ii) solves the discrete equations with high precision. The authors
are aware of two such solvers: the DOLFIN solver used in this work and
Taylor~\cite{jorba_software_2005}. The full reference
solution is available on request.

\subsection{Dual solution and stability factors}

The dual solution grows exponentially backward in time. The size of
the dual solution at time $t = 0$ is $S_D = \|z(0)\| \approx 0.510
\cdot 10^{388}$. By Theorem~\eqref{th:errorestimate}, it follows that
perturbations in initial data for the Lorenz system are amplified by a
factor $10^{388}$ at time $T = 1000$. The amplification of round-off
errors may be estimated similarly by integrating the norm of the dual
solution over the time interval. One finds that $S_C = \int_0^T \|\pi
z\| \dt \approx 2.08 \cdot 10^{388}$, which is the amplification of
errors caused by finite precision arithmetic. The stability factor for
discretization errors depends on the numerical method and in the case
of the \cg{1} method, one finds that $S_G = \int_0^T \|\dot{z}\| \dt
\approx 28.9 \cdot 10^{388}$. This is summarized in
Table~\ref{tab:stability}.

\begin{table}
  \begin{center}
    \begin{tabular}{|c|c|c|}
      \hline
      $S_D$ & $S_G$ & $S_C$ \\
      \hline
      $0.510 \cdot 10^{388}$ &  $28.9 \cdot 10^{388}$ &  $2.08 \cdot 10^{388}$ \\
      \hline
    \end{tabular}
    \caption{Size of the stability factors $S_D$, $S_G$ (for \cg{1}),
      and $S_C$ at $T = 1000$.}
    \label{tab:stability}
  \end{center}
\end{table}

\begin{figure}
  \begin{center}
    \begin{tabular}{cc}
      \begin{minipage}{5cm}
        \includegraphics[height=4cm]{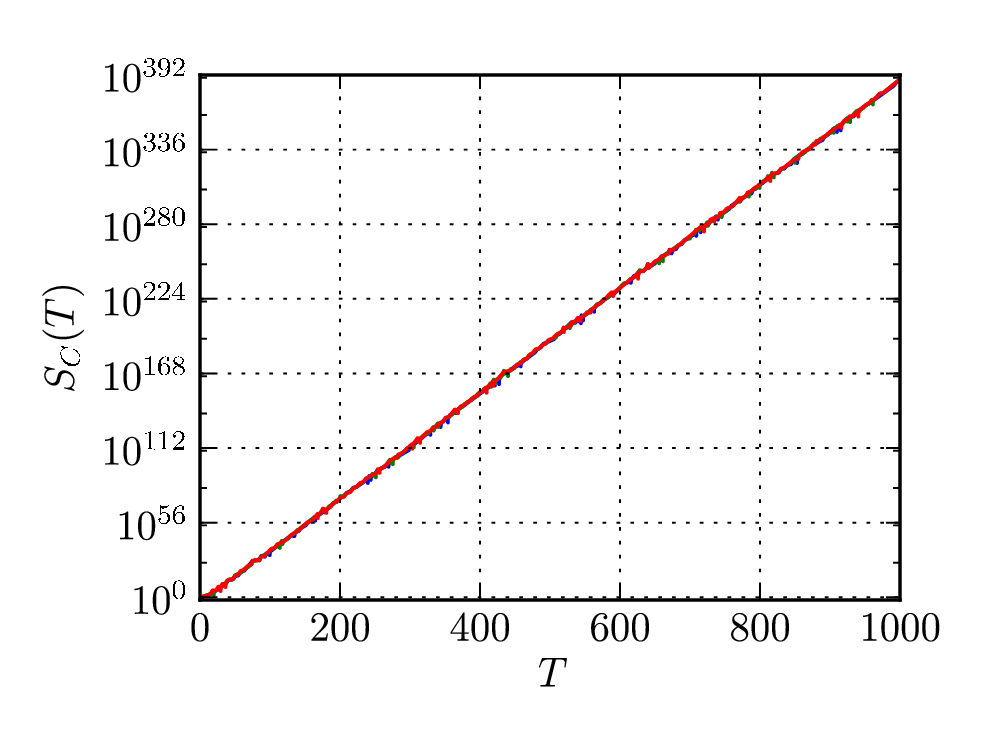}
      \end{minipage}
      &
      \begin{minipage}{5cm}
        \vspace{-0.37cm}
        \includegraphics[height=3.55cm]{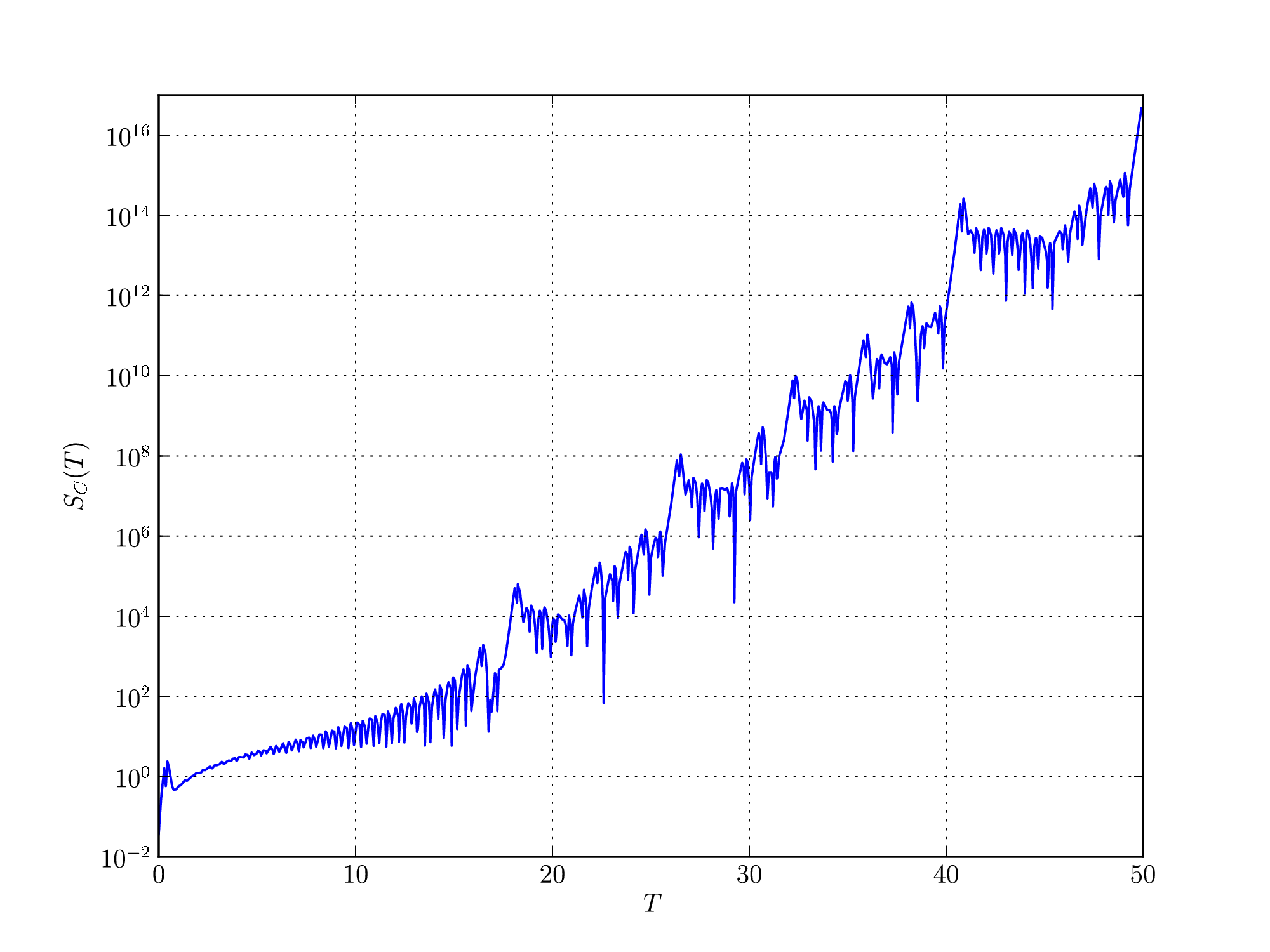}
      \end{minipage}
    \end{tabular}
    \caption{Growth of the stability factors $S_C$ on the time
      interval~$[0, 1000]$ (left) and $[0, 50]$ right.}
    \label{fig:stability}
  \end{center}
\end{figure}

By repeatedly solving the dual problem on time intervals of increasing
size, it is possible to examine the growth of the stability factors as
function of the end time~$T$. The result is displayed in
Figure~\ref{fig:stability}. Note that each data point $(T, S)$ in
Figure~\ref{fig:stability} corresponds to a solution of the dual
problem on the interval $[0, T]$.

By Figure~\ref{fig:stability}, it is evident that the stability
factors grow exponentially with the end time~$T$. On $[0, 1000]$, the
growth of the stability factor(s) may be approximated by
\begin{equation}
  S(T) \sim 10^{0.388 T} \sim 10^{0.4 T}.
\end{equation}
The rate of growth is very stable and it is therefore reasonable to
extrapolate beyond time $T = 1000$ to predict the computability of the
Lorenz system on $[0, \infty)$. We return to this question below in
  Section~\ref{sec:computability}.

A growth rate of $10^{0.388 T}$ is far below the growth rate $e^{33T}$
indicated by the simple analytic \emph{a~priori} error
estimate~\eqref{eq:apriori}. A close inspection of the growth of the
stability factor $S_C$ (Figure~\ref{fig:stability}) explains the
discrepancy between the two estimates. The growth rate of the
stability factor is not constant; it is not even monotonically
increasing. While it sometimes grows very rapidly, the average growth
rate is much smaller. The analytic \emph{a~priori} estimate must
account for the worst case growth rate and will therefore overestimate
the rate of error accumulation by a large margin.

\subsection{Error propagation}

\begin{figure}
  \begin{center}
    \includegraphics[width=6cm]{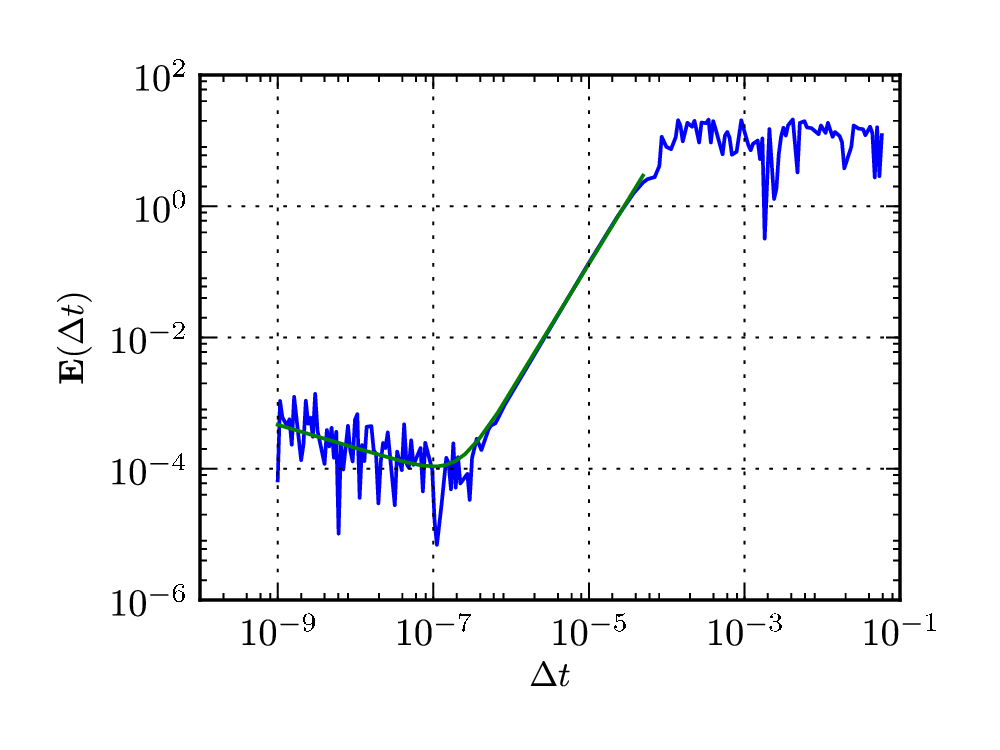}
    \includegraphics[width=6cm]{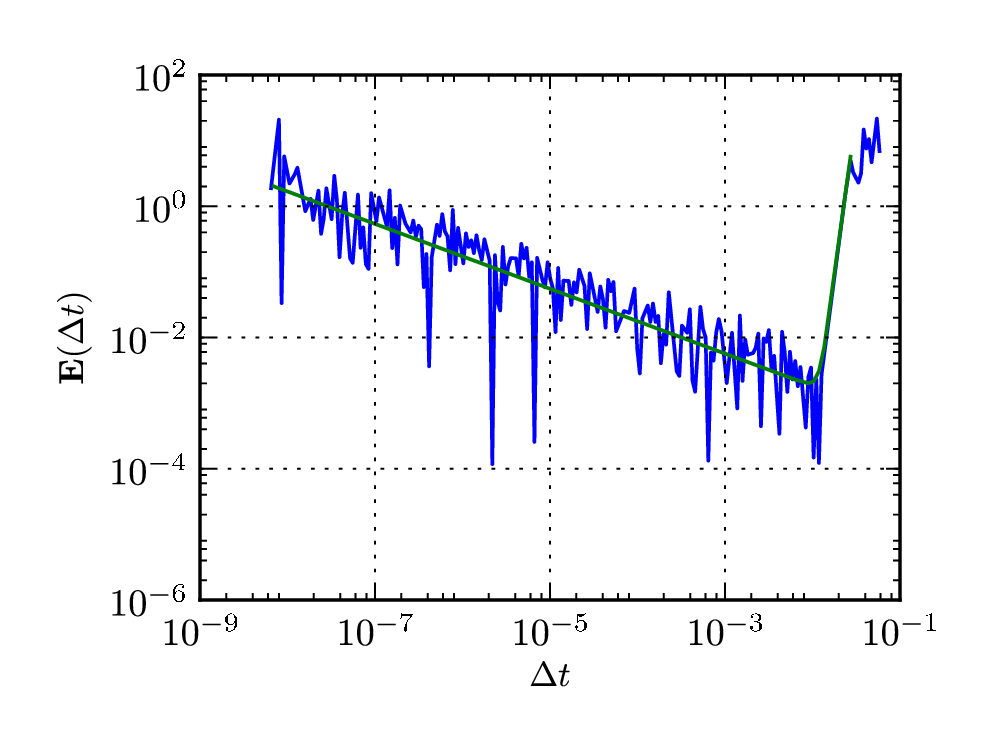}
    \caption{Error at time $T = 30$ for the \cg{1} solution (left) and
      at time $T = 40$ for the \cg{5} solution (right) of the Lorenz
      system.  The slopes of the green lines are $-0.35 \approx -1/2$
      and $1.95 \approx 2$ for the \cg{1} method. For the \cg{5}
      method, the slopes are $-0.49 \approx -1/2$ and $10.00 \approx
      10$.}
    \label{fig:errors}
  \end{center}
\end{figure}

We conclude this section by examining how the error depends on the
size of the time step $\k$. In Section~\ref{sec:erroranalysis}, we
found that the discretization error $\EG$ scales like~$\k^{2q}$ for
the \cg{q} method. On the other hand, we expect the computational
error $\EC$ to scale like~$\k^{-1/2}$. Since initial data is
represented with very high precision, we have $\EE \approx \EG + \EC
\sim \k^{2q} + \k^{-1/2}$. We thus expect the error to decrease when
the time step is decreased, at least initially. However, at the point
where $\EG = \EC$, the computational error will start to dominate and
we expect to see the error \emph{increase} with decreasing time
step. This is confirmed by the results presented in
Figure~\ref{fig:errors}, which also confirm the convergence rates $\EG
\sim \k^{2q}$ and $\EC \sim \k^{-1/2}$. We also note that the error
remains bounded for large values of $\k$; the numerical solution stays
close to the attractor but in the wrong place.

\section{Computability of the Lorenz system}
\label{sec:computability}

\subsection{A model for the computability of the Lorenz system}

Based on the analysis of Section~\ref{sec:erroranalysis} and the
numerical results of Section~\ref{sec:numericalresults}, we develop a
model for the computability of the Lorenz system. We consider the
\cg{q} method and make the following Ansatz for the error at the final
time~$T$ as function of the time step~$\k$, the polynomial degree~$q$,
and the precision~$\emach$,
\begin{equation*}
  \EE =
  \left[ C_1^{[q]} \|U(0) - u(0)\| + C_2^{[q]} \k^{\alpha} + C_3^{[q]} \k^{\beta} \emach \right]
  \cdot 10^{0.388 T}.
\end{equation*}
To determine the constants $\alpha$, $\beta$, $C_1^{[q]}$,
$C_2^{[q]}$, and $C_3^{[q]}$, we repeat the experiment of
Figure~\ref{fig:errors} for $q = 2,3,4,5$ on the interval $[0, 40]$
using the \cg{100} solution as a reference. The constants $\alpha$ and
$\beta$ may be determined by a least-squares fitting of a linear
polynomial to the regime where the error is dominated by the
discretization error $\EG$ or the computational error $\EC$,
respectively. The results are given in Table~\ref{tab:slopes}. As
expected, we find that $\alpha \approx 2q$. Furthermore, we find that
$\beta \approx -1/2$ in agreement with Theorem~\ref{th:random}.

\begin{table}
  \begin{center}
    \begin{tabular}{|c||c|c|c|c|c|c|c|c|c|c|}
      \hline
      $q$     & 2  & 3  & 4  & 5  \\
      \hline
      \hline
      $\alpha$ & 4.04 & 5.46 & 8.15 & 10.00 \\
      \hline
      $\beta$ & -0.47 & -0.50 & -0.50 & -0.49 \\
      \hline
    \end{tabular}
    \caption{Values of the constants $\alpha$ and $\beta$ as function
      of $q$ at time $T=40$.}
    \label{tab:slopes}
  \end{center}
\end{table}

Next, we fix the constants $\alpha = 2q$ and $\beta = -1/2$ and
determine the constants $C_1^{[q]}$, $C_2^{[q]}$, and $C_3^{[q]}$ as
function of $q$. In Section~\ref{sec:erroranalysis}, we found that
$S_D(T) = \|z(0)\| \approx 0.510 \cdot 10^{0.388T}$; hence $C_1^{[q]}
\approx 0.5$. By fitting curves of the form $C_2^{[q]} \k^{2q} \cdot
10^{0.388 T}$ and $C_3^{[q]} \k^{-1/2} \cdot 10^{0.388 T}$ to the two
regimes where either $\EG$ or $\EC$ dominates, we find values for the
constants $C_2^{[q]}$ and $C_3^{[q]}$. We expect $C_2^{[q]}$ to
decrease with increasing $q$ (it is essentially an interpolation
constant) and $C_3^{[q]}$ to grow at a moderate rate (by a close
inspection of the proof of Theorem~\ref{th:random}). The results are
listed in Table~\ref{tab:constants}. Based on these results, we find
that
\begin{equation*}
  \begin{split}
    C_2^{[q]} &< 0.001,  \\
    C_3^{[q]} &\approx 0.002 + 0.0005q.
  \end{split}
\end{equation*}
We thus arrive at the following model for the propagation of errors:
\begin{equation} \label{eq:model}
  \EE \approx \left[ 0.5 \, \|U(0) - u(0)\| + 0.001 \Delta t^{2q} + (0.002 + 0.0005q) \Delta t^{-1/2} \emach) \right] \cdot 10^{0.388 T}.
\end{equation}

\begin{table}
  \begin{center}
    \begin{tabular}{|c||c|c|c|c|c|c|c|c|c|c|}
      \hline
      $q$     & 2  & 3  & 4  & 5  \\
      \hline
      \hline
      $C_2^{[q]}$ & $0.000356$ & $0.000135$ & $0.000032$ & $0.000007$ \\
      \hline
      $C_3^{[q]}$ & $0.0031$ & $0.0036$ & $0.0042$ & $0.0048$ \\
      \hline
    \end{tabular}
    \caption{Values of the constants $C_2^{[q]}$ and $C_3^{[q]}$ as function
      of $q$.}
    \label{tab:constants}
  \end{center}
\end{table}

\subsection{Optimal time step}

Based on the model~\eqref{eq:model}, we determine an estimate of the
optimal time step size by setting $\EG = \EC$. We find that
\begin{equation} \label{eq:optimaltimestep}
  \k = ((2 + 0.5q) \emach)^{\frac{1}{2q + 1/2}}
  \approx
  \emach^{\frac{1}{2q + 1/2}}
\end{equation}
for large values of $q$. Inserting the values $\emach = 10^{-420}$
and $q = 100$ used in this work, we find $\k \approx 0.008$ which is
reasonably close to the value of $\k = 0.0037$ which was used to
compute the solution.

\subsection{Computability as function of machine precision}

To answer the question posed in the introduction --- \emph{How far is
  the solution computable for a given machine precision?} --- we
insert the approximate optimal time step~$\k$ given
by~\eqref{eq:optimaltimestep} into~\eqref{eq:model}. Neglecting data
errors, that is, assuming $U(0) = u(0)$, we find that
\begin{equation*}
  \EE \approx 2 \cdot 0.001 \k^{2q} \cdot 10^{0.388 T}
  \approx 0.002 \emach^{\frac{2q}{2q + 1/2}} \cdot 10^{0.388 T}
  \approx 0.002 \emach \cdot 10^{0.4 T}
\end{equation*}
for large values of $q$. Let $\nmach = -\log_{10} \emach$ be the
number of significant digits. It follows that $\EE \approx 0.002 \cdot
10^{0.4 T - \nmach}$. We conclude that the computability
$T_{\epsilon}$, that is, the time $T_{\epsilon}$ at which the solution
is no longer accurate to within a precision $\epsilon$ is
\begin{equation*}
  T_{\epsilon}(\emach) = \frac{\nmach + \log_{10}(\epsilon/0.002)}{0.4}.
\end{equation*}
Since $T_{\epsilon}$ does not depend strongly on $\epsilon$ (for
$\emach \ll \epsilon$), we find that the computability of the Lorenz
system is given by
\begin{equation*}
  T(\emach) = \nmach / 0.4 = 2.5 \, \nmach.
\end{equation*}

With six significant digits available to Lorenz in 1963, the
computability was limited to $T \approx 2.5 \cdot 6 = 15$. With
16~significant digits, the computability is limited to $T \approx 2.5
\cdot 16 = 40$. Finally, with 420~significant digits, as was used in
this work, the computability is limited to
\begin{equation*}
  T \sim 2.5 \cdot 420 = 1050 > 1000.
\end{equation*}

A more precise estimate is possible by considering the actual size of
the stability factor at any given time $T$. Noting that $S_C(T)
\approx 2 \cdot 10^{388}$ at $T = 1000$, we may obtain the estimate
\begin{equation*}
  \EE \approx 0.001 \emach \, S_C(T).
\end{equation*}
With $\emach = 10^{-16}$, it follows from Figure~\ref{fig:stability}
that $\EE = 0.001$ at $T \approx 50$. Furthermore, for $\emach =
10^{-6}$ we find that the computability is limited to $T \approx 25$.

\section{Conclusions}

We have investigated the computability of the Lorenz system and come
to the conclusion that the size of the time interval on which the
solution is computable scales linearly with the number of digits, $T
\sim 2.5 \, \nmach$. Thus, with 420 digits of precision, as used in
this work, the computability is limited to $2.5 \cdot 420 \approx
1000$.  Furthermore, if a precision of 840~digits is used, one may
compute the solution on the time interval $[0, 2000]$ and if a
precision of 4200~digits is used, one may compute the solution on the
time interval $[0, 10000]$.

\section*{Acknowledgements}
This work is supported by an Outstanding Young Investigator grant from
the Research Council of Norway, NFR 180450. This work is also
supported by a Center of Excellence grant from the Research Council of
Norway to the Center for Biomedical Computing at Simula Research
Laboratory.

\bibliography{bibliography}
\bibliographystyle{plain}


\end{document}